\documentclass[12pt]{article}
\usepackage{amsmath,amssymb,amsthm, amscd}

\begin{document}
\begin{center} 
{\bf Fundamental groups of symplectic singularities} 
\end{center} 
\vspace{0.4cm}

\begin{center}
Yoshinori Namikawa
\end{center} 
\vspace{0.4cm}

\begin{center}

{\bf Introduction}
\end{center}
 
Let $(X, \omega)$ be an affine symplectic variety. By the definition [Be], $X$ is 
an affine normal variety and $\omega$ is a holomorphic symplectic 2-form 
on the regular locus $X_{reg}$ of $X$ such that it extends to a holomorphic 
(not necessarily symplectic) 2-form on a resolution $\tilde{X}$ of $X$. 
In this article we also assume that $X$ has a $\mathbf{C}^*$-action with 
positive weights and that $\omega$ is homogeneous with respect to the 
$\mathbf{C}^*$-action. 
More precisely, the affine ring $R$ of $X$ is positively 
graded: $R = \oplus_{i \geq 0}R_i$ with $R_0 = \mathbf{C}$ and there is an 
integer $l$ such that $t^*\omega = t^l\cdot \omega$ for all $t \in \mathbf{C}^*$. 
Since $X$ has canonical singularities, we have $l > 0$ ([Na, Lemma (2.2)]).
Affine symplectic varieties are constructed in various ways such as nilpotent orbit 
closures of a semisimple complex Lie algebra (cf. [C-M]), Slodowy slices to nilpotent 
orbits ([Sl]) or symplectic reductions of holomorphic symplectic manifolds with 
Hamiltonian actions. These varieties come up with $\mathbf{C}^*$-actions and   
the above assumption of the $\mathbf{C}^*$-action is satisfied in all examples we know. 
  
In the previous article [Na] we posed a question: 
\vspace{0.2cm}

{\bf Problem}. {\em Is the fundamental group $\pi_1(X_{reg})$ finite ?} 
\vspace{0.2cm}

Such fundamental groups are explicitly calculated by a group-theoretic method 
when $X$ is a nilpotent orbit closure (cf. [C-M]). However no general results are 
known.  
 
In this short note we give a partial answer to this question. Namely we have 
\vspace{0.2cm}

{\bf Theorem}. {\em The algebraic fundamental group $\hat{\pi}_1(X_{reg})$ is a finite 
group.} 
\vspace{0.2cm}  

Notice that a symplectic variety $X$ has canonical singularities. In particular, the 
log pair $(X, 0)$ has klt (Kawamata log terminal) singularities.  
The theorem is, in fact, a corollary to the more general result:
\vspace{0.2cm}

{\bf Main Theorem}. {\em Let $X := \mathrm{Spec} R$ be an affine variety where 
$R$ is positively graded:  $R = \oplus_{i \geq 0}R_i$ with $R_0 = \mathbf{C}$. 
Assume that the log pair $(X, 0)$ has klt singulatities. 
Then $\hat{\pi}_1(X_{reg})$ is a finite group.}      
\vspace{0.2cm}

Recently C. Xu [Xu] proved that, for a klt pair $(X, \Delta)$ and a point $p \in X$, the 
algebraic fundamental group $\hat{\pi}_1(U - \{p\})$ is finite for a small complex analytic 
neighborhood $U$ of $p$. By using the Koll\'{a}r component he obtained it 
from the finiteness of the algebraic fundamental group of the regular part of a 
log Fano variety. 
Since $X$ has a $\mathbf{C}^*$-action in our case, 
$\hat{\pi}_1(X_{reg}) \cong \hat{\pi}_1(U_{reg})$ 
for a small complex analytic neighborhood $U$ of the origin 
$p \in X$. The argument in [Xu] is also valid for $\hat{\pi}_1(U_{reg})$ and one can 
prove Main Theorem. 
 
In this article we introduce another approach to Main Theorem by using the orbifold 
fundamental group.  
 
To explain the basic idea of the proof we first assume that $R$ is generated by $R_1$ as 
a $\mathbf{C}$-algebra and $X$ has only isolated singularity. Put $\mathbf{P}(X) := 
\mathrm{Proj}R$. By the assumption $\mathbf{P}(X)$ is a projective manifold. Since 
$(X,0)$ has klt singularities, we also see that $\mathbf{P}(X)$ is a Fano manifold (cf. 
[Fu], Prposition 4.38). 
Let $L$ be the tautological line bundle on $\mathbf{P}(X)$ and 
denote by $(L^{-1})^{\times}$ be the $\mathbf{C}^*$-bundle on $\mathbf{P}(X)$ 
obtained from $L^{-1}$ by removing the 0-section. Then the projection map 
$p: X - \{0\} \to \mathbf{P}(X)$ can be identified with $(L^{-1})^{\times} \to \mathbf{P}(X)$. 
There is a homotopy exact sequence 
$$ \pi_1(\mathbf{C}^*) \to \pi_1(X - \{0\}) \to \pi_1(\mathbf{P}(X)) 
\to 1. $$ Here $\pi_1(\mathbf{P}(X)) = 1$ 
because $\mathbf{P}(X)$ is a Fano manifold. We want to show that  
the first map $\pi_1(\mathbf{C}^*) \to \pi_1(X - \{0\})$ has a nontrivial kernel. 
Suppose to the contrary that it is an injection. Then $\pi_1(X - \{0\}) = \mathbf{Z}$ 
and one has a surjective map $\pi_1(X - \{0\}) \to \mathbf{Z}/l\mathbf{Z}$ for 
any $l > 1$. This determines an \'{e}tale covering $f: Y \to X - \{0\}$, which extends to 
a finite surjective map $\bar{f}: \bar{Y} \to L^{-1}$, where $\bar{Y}$ contains $Y$ as 
a Zariski open subset and $\bar{f}$ is a cyclic covering branched along the 0-section $\Sigma$  
of $L^{-1}$. The direct image $\bar{f}_*O_{\bar Y}$ can be written as $O_{\bar Y} 
\oplus M \oplus M^{\otimes 2} \oplus ... \oplus M^{\otimes l -1}$ with a line bundle $M$ 
on $L^{-1}$. Here $M^{\otimes l} \cong O_{L^{-1}}(-\Sigma)$. Restrict this isomorphism 
to $\Sigma (\cong \mathbf{P}(X))$. Then we have $(M\vert_{\Sigma})^{\otimes l} \cong L$. 
This shows that $L \in \mathrm{Pic}(\mathbf{P}(X))$ is divisible by any $l > 1$. 
But this is absurd because $L$ is an ample line bundle.    
Therefore $\pi_1(X - \{0\})$ is finite. 

In a general situation $\mathbf{P}(X)$ is no more smooth and the projection map 
$X - \{0\} \to \mathbf{P}(X)$ is not a $\mathbf{C}^*$-bundle. We take a smooth open set 
$\mathbf{P}(X)^{\sharp}$ of $\mathbf{P}(X)$ in such a way that 
$X^{\sharp} := p^{-1}(\mathbf{P}(X)^{\sharp})$ is smooth and 
$\mathrm{Codim}_{\mathbf{P}(X)}(\mathbf{P}(X) - \mathbf{P}(X)^{\sharp}) \geq 2$. 
The map $X^{\sharp} \to \mathbf{P}(X)^{\sharp}$ is not still a $\mathbf{C}^*$-bundle, but 
if we introduce a suitable orbifold structure on $\mathbf{P}(X)$, then it can 
be regarded as a $\mathbf{C}^*$-bundle on the orbifold $\mathbf{P}(X)^{\sharp, orb}$. 
Moreover we have a homotopy exact sequence 
$$ \pi_1(\mathbf{C}^*) \to 
\pi_1(X^{\sharp}) \to \pi_1^{orb}(\mathbf{P}(X)^{\sharp, orb}) \to 1.$$ 
The orbifold structure on $\mathbf{P}(X)$ determines an effective $\mathbf{Q}$-divisor 
$\Delta$ with standard coefficients. By the assumption that $(X, 0)$ has klt singularities, 
we see that $(\mathbf{P}(X), \Delta)$ is a log Fano variety (\S 1, Lemma). C. Xu [Xu] has proved that $\hat{\pi}_1(\mathbf{P}(X)_{reg})$ is finite for such a variety. It turns out that 
his proof can be used to prove that $\hat{\pi}_1^{orb}(\mathbf{P}(X)^{\sharp, orb})$ 
is finite.  We take a finite \'{e}tale covering $Y^{orb} \to \mathbf{P}(X)^{\sharp, orb}$ 
such that $\hat{\pi}_1^{orb}(Y^{orb}) = 1$ and define $Z$ to be the normalization of 
$X^{\sharp} \times_{\mathbf{P}(X)^{\sharp}}Y$. Then we see that $Z \to X^{\sharp}$ is 
an \'{e}tale covering in the usual sense. Moreover, we have an exact 
sequence 
$$ \pi_1(\mathbf{C}^*) \to 
\pi_1(Z) \to \pi_1^{orb}(Y^{orb}) \to 1$$ by replacing $X^{\sharp}$ and $\mathbf{P}(X)^{\sharp, orb}$ by $Z$ and $Y^{orb}$. Assume that there exists a surjection from $\pi_1(Z)$ to a finite 
group $\Gamma$. Since $\pi_1^{orb}(Y^{orb})$ has no nontrivial finite quotient, the 
composition map $\pi_1(\mathbf{C}^*) \to \pi_1(Z) \to \Gamma$ is surjective.  
The orbifold line bundle associated with the orbifold $\mathbf{C}^*$-bundle $Z \to Y^{orb}$ is 
negative. We prove that the order of $\Gamma$ cannot be arbitrary 
large by using this fact; hence $\hat{\pi}_1(Z)$ is finite. Since $\hat{\pi}_1(Z)$ is a finite index subgroup of $\hat{\pi}_1(X^{\sharp})$, $\hat{\pi}_1(X^{\sharp})$ is also finite.  
As there is a surjection map $\hat{\pi}_1(X^{\sharp}) \to \hat{\pi}_1(X_{reg})$, we have Main  Theorem. 

The argument above also shows that $\pi_1(X_{reg})$ is finite if and only if 
$\pi_1^{orb}(\mathbf{P}(X)^{\sharp, orb})$ is finite.  
    
\vspace{0.2cm}   

{\em Acknowledgement}: The author thanks Y. Kawamata and Y.Gongyo for pointing 
out that the klt condition would be enough for proving Theorem. 
\vspace{0.2cm}

{\bf 1. Algebraic orbifolds} 

In the remainder of this article $X := \mathrm{Spec} R$ is an affine normal variety with a  
positively graded ring $R = \oplus_{i \geq 0}R_i$, $R_0 = \mathbf{C}$ such that $(X,0)$ is 
a klt pair.
Take minimal homogeneous generators of $R$ and consider  
the surjection 
$$\mathbf{C}[x_0, ..., x_n] \to R$$ which sends each $x_i$ to the homogeneous 
generator.  Correspondingly $X$ is embedded in 
$\mathbf{C}^{n+1}$. To $x_i$ we give the same weight as the minimal generator. 
Put  $a_i := wt(x_i)$. We may assume that $GCD(a_0, ..., a_n) = 1$.  
The quotient variety 
$\mathbf{C}^{n+1}-\{0\}/\mathbf{C}^*$ by the $\mathbf{C}^*$-action 
$(x_0, ..., x_n) \to (t^{a_0}x_0, ..., t^{a_n}x_n)$ is the weighted projective space 
$\mathbf{P}(a_0, ..., a_n)$. We put $\mathbf{P}(X) := X - \{0\}/\mathbf{C}^*$. 
By the definition $\mathbf{P}(X)$ is a closed subvariety of 
$\mathbf{P}(a_0, ..., a_n)$.  Put $W_i := \{ x_i = 1\} \subset \mathbf{C}^{n+1}$. 
Then the projection map $p: \mathbf{C}^{n+1}- \{0\} \to \mathbf{P}(a_0, ..., a_n)$ 
induces a map $p_i: W_i \to \mathbf{P}(a_0, ..., a_n)$, which is a finite Galois covering 
of the image. The collection $\{p_i\}$ defines a smooth orbifold structure on 
$\mathbf{P}(a_0, ..., a_n)$ in the sense of [Mum, \S 2]. 
More exactly, the following are satisfied 

(i) For each $i$, $W_i$ is a smooth variety and $p_i: W_i \to p_i(W_i)$ 
is a finite Galois covering\footnote{The precise definition of an orbifold only needs a 
slightly weaker condition: $p_i: W_i \to \mathbf{P}(a_0, ..., a_n)$ factorizes as 
$W_i \stackrel{q_i}\to W_i/G_i \stackrel{r_i}\to \mathbf{P}(a_0, ..., a_n)$ where $G_i$ is a 
finite group and $r_i$ is an \'{e}tale map.}.  
$\cup \mathrm{Im}(p_i) = \mathbf{P}(a_0, ..., a_n)$.    

(ii) Let $(W_i \times_{\mathbf{P}(a_0, ..., a_n)}W_j)^n$ denote the normalization 
of the fibre product $W_i \times_{\mathbf{P}(a_0, ..., a_n)}W_j$. Then 
the maps $(W_i \times_{\mathbf{P}(a_0, ..., a_n)}W_j)^n \to W_i$ and 
$(W_i \times_{\mathbf{P}(a_0, ..., a_n)}W_j)^n \to W_j$ are both \'{e}tale maps.     
\vspace{0.2cm}

The orbifold $\mathbf{P}(a_0, ..., a_n)$ 
admits an orbifold line bundle $O_{\mathbf{P}(a_0, ..., a_n)}(1)$.  
Put $D_i := \{x_i = 0\} \subset \mathbf{P}(a_0, ..., a_n)$ and 
$D := \cup D_i$. Since $x_i$ are minimal generators, $\bar{D}:= \mathbf{P}(X) \cap D$ 
is a divisor of $\mathbf{P}(X)$. Define 
$$\mathbf{P}(X)^{\sharp} := \mathbf{P}(X) - \mathrm{Sing}(\bar{D}) - \mathrm{Sing}(\mathbf{P}(X)),$$ and  
$$X^{\sharp} := p^{-1}(\mathbf{P}(X)^{\sharp}).$$
 
Let $$\bar{D}= \cup \bar{D}_{\alpha}$$ 
be the decomposition into irreducible components \footnote{The index $\alpha$ 
is usually different from the original index $i$ of $D_i$ because $D_{i_1} \cap ... 
\cap D_{i_k} \cap \mathbf{P}(X)$ may possibly become an irreducible component 
of $\bar{D}$ or $D_i \cap \mathbf{P}(X)$ may split into more than two irreducible 
components of $\bar{D}$.}.   
By the definition $\bar{D}^{\sharp} := \bar{D} \cap \mathbf{P}(X)^{\sharp}$ 
is a smooth divisor of $\mathbf{P}(X)^{\sharp}$.
Put $\bar{D}^{\sharp}_{\alpha} := \bar{D}_{\alpha} \cap \mathbf{P}(X)^{\sharp}$. 
Then $\bar{D}^{\sharp}$ is the disjoint union of irreducible smooth divisors 
$\bar{D}^{\sharp}_{\alpha}$.  
    
In general $p^{-1}(\mathbf{P}(X)_{reg})$ is not smooth; but if we shrink 
$\mathbf{P}(X)_{reg}$ to $\mathbf{P}(X)^{\sharp}$, then its inverse image  
$X^{\sharp}$ is smooth. 

Notice that every fibre of $p^{\sharp} (:= p\vert_{X^{\sharp}}) : X^{\sharp} 
\to \mathbf{P}(X)^{\sharp}$ is isomorphic to $\mathbf{C}^*$, but the fibre  
over a point of $\bar{D}^{\sharp}_{\alpha}$ may possibly be a multiple fibre. 
We denote by $m_{\alpha}$\footnote{The multiplicity $m_{\alpha}$ may possibly be one.} the multiplicity of a fibre over a point of 
$\bar{D}^{\sharp}_{\alpha}$. 
  
The map $p^{\sharp}$ is a $\mathbf{C}^*$-bundle if we restrict it to the 
open set $\mathbf{P}(X)^{\sharp} - \bar{D}^{\sharp}$. 
Notice that  $\mathbf{P}(X) - \mathbf{P}(X)^{\sharp}$ has at least codimension 2 in 
$\mathbf{P}(X)$.   

By putting $U_i := X \cap W_i$ and $\pi_i := p_i\vert_{U_i}$, the collection  
$\{\pi_i: U_i \to \mathbf{P}(X)\}$ of covering maps induces a (not necessarily smooth) 
orbifold structure on $\mathbf{P}(X)$.    
Namely, we have 

(i) For each $i$, $U_i$ is a normal variety and $\pi_i: U_i \to \pi_i(U_i)$ is a finite 
Galois covering. $\cup \mathrm{Im}(\pi_i) = \mathbf{P}(X)$  

(ii) The maps $(U_i \times_{\mathbf{P}(X)}U_j)^n \to U_i$ and 
$(U_i \times_{\mathbf{P}(X)}U_j)^n \to U_j$ are both \'{e}tale maps.
\vspace{0.2cm}

We put 
$\mathcal{L} := O_{\mathbf{P}(a_0, ..., a_n)}(1)\vert_{\mathbf{P}(X)}$, which is an 
orbifold line bundle on $\mathbf{P}(X)$. We call $\mathcal{L}$ the tautological 
line bundle. Then $X - \{0\} \to \mathbf{P}(X)$ can be regarded as an orbifold line 
bundlle $\mathcal{L}^{-1}$  

Notice that, if we restrict this orbifold structure to $\mathbf{P}(X)^{\sharp}$, then 
it is a smooth orbifold structure.  
 \vspace{0.2cm}

{\bf Lemma}. {\em Assume that the log pair $(X,0)$ has klt singularities. Put 
$\Delta := \sum (1 - 1/m_{\alpha})\bar{D}_{\alpha}$. Then 
$(\mathbf{P}(X), \Delta)$ is a log Fano variety, that is, $(\mathbf{P}(X), \Delta)$ 
has klt singularities and $-(K_{\mathbf{P}(X)} + \Delta)$ is an ample $\mathbf{Q}$-divisor.} 
\vspace{0.2cm}

{\em Proof}.  Take a positive integer $d$ in such a way that  
the subring $R^{(d)} := \oplus_{i \geq 0}R_{id}$ is generated by $R_d$ as a $\mathbf{C}$-algebra. We put $V := \mathrm{Spec} R^{(d)}$. Then there is a finite 
surjective map $\mu: X \to V$. Notice that $\mathrm{Proj}(R) = \mathrm{Proj}(R^{(d)})$. 
Hence there is a natural projection map $q: V - \{0\} \to \mathbf{P}(X)$ and 
the composition map $X - \{0\} \to V - \{0\} \to \mathbf{P}(X)$ coincides with the natural 
projection map $p: X - \{0\} \to \mathbf{P}(X)$.  
Since $R^{(d)}$ is generated by $R_d$ as a $\mathbf{C}$-algebra, the projection map $q$ is a $\mathbf{C}^*$-bundle. We put $V^{\sharp} := q^{-1}(\mathbf{P}(X)^{\sharp})$ and put $q^{\sharp} := q\vert_{V^{\sharp}}: V^{\sharp} \to \mathbf{P}(X)^{\sharp}$.
For a point $t \in \mathbf{P}(X)^{\sharp} - \bar{D}^{\sharp}$, the fibres 
$(p^{\sharp})^{-1}(t)$ and $(q^{\sharp})^{-1}(t)$ are both isomorphic to $\mathbf{C}^*$ and $\mu$ induces an 
etale covering between them of the same degree as $\mathrm{deg}(\mu)$. 
On the other hand, for a point $t \in \bar{D}_{\alpha}^{\sharp}$, the fibre $(p^{\sharp})^{-1}(t)$ is a multiple fibre with multiplicity $m_{\alpha}$ and $(p^{\sharp})^{-1}(t)_{red} \cong \mathbf{C}^*$. 
In this case $\mu$ induces an etale covering  $(p^{\sharp})^{-1}(t)_{red} \to (q^{\sharp})^{-1}(t)$ of degree $\mathrm{deg}(\mu)/m_{\alpha}$. In other words, $X^{\sharp} \to V^{\sharp}$ is 
a finite cover, which is branched along $(q^{\sharp})^{-1}(\cup_{\alpha; m_{\alpha} > 1}\bar{D}^{\sharp}_{\alpha})$. 
Let $B$ be the $\mathbf{Q}$-divisor of $V$ obtained as the closure of the $\mathbf{Q}$-divisor $q^*\Delta$ of $V - \{0\}$. 
Here  notice that $\mathrm{Codim}_X(X - X^{\sharp}) \geq 2$ and 
$\mathrm{Codim}_V(V - V^{\sharp}) \geq 2$. Then we have  
$$ K_X = \mu^*(K_V + B).$$ 
Since $(X, 0)$ is a klt pair, $(V, B)$ is a klt pair by [K-M, Proposition 5.20]. 
By [Fu, Proposition 4.38] we conclude that $(\mathbf{P}(X), \Delta)$ is a log Fano variety. 
\vspace{0.2cm}

{\bf Remark}. The lemma is rephrased as:  if $(X, 0)$ has klt singularities, then 
$\mathbf{P}(X)^{orb}$ is a Fano orbifold.  
When $X$ is an affine symplectic variety, this can be proved directly by using the 
fact that $\mathbf{P}(X)$ has a contact orbifold structure ([Na, Theorem 4.4.1]). 
\vspace{0.2cm}  

{\bf 2. Algebraic orbifold fundamental group of $\mathbf{P}(X)^{\sharp}$} 
\vspace{0.2cm}

In the previous section we observed that  $\mathbf{P}(X)^{\sharp}$ has a 
smooth orbifold structure.  
Namely, if we put $U_i^{\sharp} := \pi_i^{-1}(\mathbf{P}(X)^{\sharp})$, 
and $\pi_i^{\sharp} := \pi_i\vert_{U_i^{\sharp}}$, then 
$\mathcal{U} := \{\pi_i^{\sharp}: U_i^{\sharp} \to \mathbf{P}(X)^{\sharp}\}_{i \in I}$ give the orbifold charts of $\mathbf{P}(X)^{\sharp}$. 
Assume that $\mathbf{P}(X)^{\sharp}$ has another orbifold charts 
$\mathcal{U}' := \{\pi'_j: U'_j  \to \mathbf{P}(X)^{\sharp}\}_{j \in J}$. 
Then $\mathcal{U}$ and 
$\mathcal{U}'$ are equivalent if, for each $i \in I$ and $j \in J$, two maps 
$U_i^{\sharp} \to \mathbf{P}(X)^{\sharp}$ and $U'_j \to \mathbf{P}(X)^{\sharp}$ are {\em admissible} to each other:  
in other words,     
$(U_i^{\sharp} \times_{\mathbf{P}(X)^{\sharp}}U'_j)^n \to U_i^{\sharp}$ and 
$(U_i^{\sharp} \times_{\mathbf{P}(X)^{\sharp}}U'_j)^n \to U'_j$ are both \'{e}tale maps. 
An orbifold structure on $\mathbf{P}(X)^{\sharp}$ is precisely 
an equivalence class of orbifold charts of $\mathbf{P}(X)^{\sharp}$.    
In the remainder we will denote by $\mathbf{P}(X)^{\sharp, orb}$ the orbifold structure 
defined in the previous section. 
Let $Y^{orb}$ be a smooth orbifold; namely it is a pair of a normal algebraic variety $Y$ 
and an equivalence class of orbifold charts $\mathcal{V} = \{\nu_k: V_k \to Y\}_{k \in K}$. 
Let  $f: Y \to \mathbf{P}(X)^{\sharp}$ be a finite surjective morphism of algebraic varieties. 
We say that $f$ is an \'{e}tale covering map from $Y^{orb}$ to $\mathbf{P}(X)^{\sharp, orb}$ if the 
following property holds: 

For any $k \in K$ and $i \in I$, two maps $f \circ \nu_k: V_k  \to \mathbf{P}(X)^{\sharp}$ and $\pi_i: U_i^{\sharp} \to \mathbf{P}(X)^{\sharp}$ are admissible to each other.  

Notice that $f$ is not necessarily an \'{e}tale covering map in the usual sense even if 
$f$ is an \'{e}tale covering map of orbifolds. An \'{e}tale covering map is said to be Galois 
if the underlying morphism is Galois in the usual sense. 
\vspace{0.2cm}

{\bf Lemma}. {\em For any finite \'{e}tale covering $f: Y^{orb} \to \mathbf{P}(X)^{\sharp, orb}$, 
there exists an \'{e}tale Galois covering $g: Z^{orb} \to \mathbf{P}(X)^{\sharp, orb}$ 
such that $g$ factorizes as $Z^{orb} \to Y^{orb} \stackrel{f}\to \mathbf{P}(X)^{\sharp, orb}$.} 
\vspace{0.1cm}

{\em Proof}. Let $K$ and $L$ be the function fields of $\mathbf{P}(X)^{\sharp}$ and 
$Y$. Let $M$ be the Galois closure of $L/K$ and take the normalization $Z$ of $Y$ in 
$M$. We shall give an orbifold structure on $Z$ in such a way that $Z^{orb} \to \mathbf{P}(X)^{\sharp, orb}$ is an \'{e}tale covering and it factorizes through $Y^{orb}$. 
Let $V \to Y$ be an orbifold chart of $Y^{orb}$. By the definition $V$ is smooth and 
the composition map $V \to Y \to \mathbf{P}(X)^{\sharp}$ is admissible for any orbifold chart $U \to \mathbf{P}(X)^{\sharp}$. In other words, $(V \times_{\mathbf{P}(X)^{\sharp}}U)^n 
\to V$ and $(V \times_{\mathbf{P}(X)^{\sharp}}U)^n \to U$ are both \'{e}tale maps. 
Let $L'$ be the function field of $V$ and let $M'$ be the Galois closure of $L'/K$. 
Let $W$ be the normalization of $V$ in $M'$. Then, since $M \subset M'$, there is a 
natural map $q_V: W \to Z$. It is clear that 
$\cup \mathrm{Im}(q_V) = Z$ when  $V$ runs through all orbifold charts of $Y^{orb}$.  
We prove that the map $W \to V$ is \'{e}tale. 

Here we recall an explicit construction of $W$. Assume that 
$V \to \mathbf{P}(X)^{\sharp}$ is not Galois. Then there is an element $\sigma \in 
G := \mathrm{Gal}(M'/K)$ such that $(L')^{\sigma} \ne L'$, where $(L')^{\sigma} := 
\sigma (L')$. 
We take an irreducible 
component $V_1$ of $(V \times_{\mathbf{P}(X)^{\sharp}}V^{\sigma})^n$ in such a way that 
the induced map $V_1 \to V$ is a finite covering with deg $\ge 2$.  
Let $L_1$ be the function field of $V_1$ and if $L_1/K$ is not still a Galois extension, 
we take the Galois closure $M_1$ of $L_1/K$. There exists an element $\sigma_1 
\in G_1 := \mathrm{Gal}(M_1/K)$ such that $(L_1)^{\sigma_1} \ne L_1$. 
We take an irreducible component $V_2$ of $(V_1 \times_{\mathbf{P}(X)^{\sharp}} (V_1)^{\sigma_1})^n$ in such a way that the induced map $V_2 \to V_1$ has degree 
$\geq 2$. When we repeat this process, we finally reach the $W$. 

Thus, to prove that $W$ is \'{e}tale over $V$, we only have to show that 
$(V \times_{\mathbf{P}(X)^{\sharp}} V^{\sigma})^n \to V$ and         
$(V \times_{\mathbf{P}(X)^{\sharp}} V^{\sigma})^n \to V^{\sigma}$ are 
both \'{e}tale maps. In fact, if this is proved, then $V_1 \to V$ is an \'{e}tale map and 
hence $V_1 \to \mathbf{P}(X)^{\sharp}$ and $U \to \mathbf{P}(X)^{\sharp}$ are 
admissible to each other. We may then replace $V$ by $V_1$ and continue.  
 
Before starting the proof we notice that $\sigma$ induces a $\mathbf{P}(X)^{\sharp}$-isomorphism $U \cong U^{\sigma}$. 
Since $V \to \mathbf{P}(X)^{\sharp}$ and $U \to \mathbf{P}(X)^{\sharp}$ are 
admissible to each other, $(V \times_{\mathbf{P}(X)^{\sharp}} U)^n \to U$ and 
$(V \times_{\mathbf{P}(X)^{\sharp}} U)^n \to V$ are both  
\'{e}tale maps. Since            
$V^{\sigma} \to \mathbf{P}(X)^{\sharp}$ and $U^{\sigma} \to \mathbf{P}(X)^{\sharp}$ 
are also admissible to each other, 
$(V^{\sigma} \times_{\mathbf{P}(X)^{\sharp}} U^{\sigma})^n \to U^{\sigma}$ 
and $(V^{\sigma} \times_{\mathbf{P}(X)^{\sharp}} U^{\sigma})^n \to V^{\sigma}$
are \'{e}tale maps. Here, identifying $U^{\sigma}$ with $U$ by the above isomorphism, 
we get two maps $(V^{\sigma} \times_{\mathbf{P}(X)^{\sharp}} U)^n \to U$ and 
$(V^{\sigma} \times_{\mathbf{P}(X)^{\sharp}} U)^n \to V^{\sigma}$.  
Now we have a commutative diagram 

\begin{equation} 
\begin{CD} 
(V \times_{\mathbf{P}(X)^{\sharp}}U \times_{\mathbf{P}(X)^{\sharp}} V^{\sigma})^n 
 @>>> (V^{\sigma} \times_{\mathbf{P}(X)^{\sharp}} U)^n \\ 
@VVV @VVV \\ 
(V \times_{\mathbf{P}(X)^{\sharp}}U)^n @>>> U,       
\end{CD} 
\end{equation}
where all maps are \'{e}tale. 
Let us consider the map 
$(V \times_{\mathbf{P}(X)^{\sharp}}U \times_{\mathbf{P}(X)^{\sharp}} V^{\sigma})^n 
\to V$. Since  this map factorizes through $(V \times_{\mathbf{P}(X)^{\sharp}}U)^n$, 
it is an \'{e}tale map. On the other hand, this map also factorizes as 
$$(V \times_{\mathbf{P}(X)^{\sharp}}U \times_{\mathbf{P}(X)^{\sharp}} V^{\sigma})^n 
\to (V \times_{\mathbf{P}(X)^{\sharp}} V^{\sigma})^n \to V.$$ 
As the first map is \'{e}tale, the second map 
$(V \times_{\mathbf{P}(X)^{\sharp}} V^{\sigma})^n \to V$ is an \'{e}tale map. 
By a similar reasoning we see that  
$(V \times_{\mathbf{P}(X)^{\sharp}} V^{\sigma})^n \to V^{\sigma}$ is an \'{e}tale 
map.  

We finally prove that $\{q_V: W \to Z\}$ gives an orbifold structure to $Z$. 
Since $M'$ is a Galois extension of $M$, we see that $q_V(W) $ is the quotient 
variety of $W$ by $\mathrm{Gal}(M'/M)$. Moreover, $W$ is a smooth variety 
because it is an \'{e}tale cover of a smooth variety $V$. In the remainder we 
shall check that $q_{V'}: W' \to Z$ and $q_V: W \to Z$ are admissible to each 
other. We have a commutative diagram 

\begin{equation} 
\begin{CD} 
(W \times_Y W')^n 
 @>>> W\\ 
@VVV @VVV \\ 
(V \times_Y V')^n @>>> V.       
\end{CD} 
\end{equation}
Here two vertical maps are \'{e}tale because $W \to V$ and $W' \to V'$ 
are \'{e}tale, and the second horizontal map is also \'{e}tale because 
$V \to Y$ and $V' \to Y$ are admissible to each other. 
Hence the first horizontal map $(W \times_Y W')^n \to W$ is \'{e}tale 
by the commutative diagram. The map $(W \times_Z W')^n \to W$ is 
factorizes as 
$$(W \times_Z W')^n \to (W \times_Y W')^n \to W.$$ 
Since first map is an open immersion, we see that 
$(W \times_Z W')^n \to W$ is an \'{e}tale map. Similarly, 
$(W \times_Z W')^n \to W'$ is an \'{e}tale map.  Q.E.D.   
\vspace{0.2cm}
 
Take a point $x \in \mathbf{P}(X)^{\sharp}$ in such a way that 
$x \notin \bar{D}^{\sharp}$. Then $f^{-1}(x)$ consists of exactly $\mathrm{deg}(f)$ points. 
Consider all pairs $(Y^{orb}, y)$ of \'{e}tale coverings $Y^{orb}$ 
of $\mathbf{P}(X)^{\sharp, orb}$ and $y \in Y$ lying on $x \in 
\mathbf{P}(X)^{\sharp}$. A morphism $h: (Z^{orb}, z) \to (Y^{orb},y)$ is a $\mathbf{P}(X)^{\sharp}$-morphism 
$h: Z \to Y$ with $h(z) = y$  such that it is an \'{e}tale covering map from $Z^{orb}$ to $Y^{orb}$. 
When $Z$ and $Y$ are both Galois coverings of $\mathbf{P}(X)^{\sharp}$, $h$ induces a surjective map   
$\mathrm{Aut}(Z/\mathbf{P}(X)^{\sharp}) \to \mathrm{Aut}(Y/\mathbf{P}(X)^{\sharp})$.  
As in the usual situation, we can define the algebraic orbifold fundamental group $\hat{\pi}_1^{orb}(\mathbf{P}(X)^{\sharp, orb}, x)$ 
as the profinite group $\lim \mathrm{Aut}(Y/\mathbf{P}(X)^{\sharp})$,  
where $f$ runs through all finite \'{e}tale Galois coverings of $\mathbf{P}(X)^{\sharp, orb}$.   
\vspace{0.2cm}

{\bf Theorem}([Xu]). {\em $\hat{\pi}_1^{orb}(\mathbf{P}(X)^{\sharp, orb}, x)$ is a finite 
group.} 
\vspace{0.2cm}

{\em Proof}. Write $\Delta^{\sharp}$ for $\Delta\vert_{\mathbf{P}(X)^{\sharp}}$.   
Let $Y^{orb}$ be a finite \'{e}tale covering map of $\mathbf{P}(X)^{\sharp, orb}$ and 
let $f: Y \to \mathbf{P}(X)^{\sharp}$ be the underlying map. 
We shall prove that $K_Y + \Delta_Y = f^*(K_{\mathbf{P}(X)^{\sharp}} + \Delta^{\sharp})$ 
for some effective divisor $\Delta_Y$ on $Y$. 
Let $\mathcal{V} = \{\nu_k: V_k 
\to Y\}_{k \in K}$ be orbifold covering charts. Let $Z$ be an irreducible component of the normalization $(V_k \times_{\mathbf{P}(X)^{\sharp}}U_i^{\sharp})^n$ of the 
fibre product of the diagram $$ V_k \stackrel{f \circ \nu_k}\to \mathbf{P}(X)^{\sharp} \leftarrow U_i^{\sharp}. $$  Then we have a commutative diagram 
\begin{equation} 
\begin{CD} 
Z @>{p_2}>> U_i^{\sharp} \\ 
@V{p_1}VV @V{\pi_i}VV \\ 
V_k @>{f \circ \nu_k}>> \mathbf{P}(X)^{\sharp}.      
\end{CD} 
\end{equation}

Here $p_1$ and $p_2$ are both \'{e}tale maps. 
Since $K_{U_i^{\sharp}} = \pi_i^*(K_{\mathbf{P}(X)^{\sharp}} + \Delta^{\sharp})$ 
and $K_Z = p_2^*K_{U_i^{\sharp}}$, we have 
$K_Z = (\pi_i \circ p_2)^*(K_{\mathbf{P}(X)^{\sharp}} + \Delta^{\sharp})$. 
On the other hand, since $K_Z = p_1^*K_{V_k}$, we see that 
$K_{V_k} = (f \circ \nu_k)^*(K_{\mathbf{P}(X)^{\sharp}} + \Delta^{\sharp})$. 
Then one can write $K_Y + \Delta_Y = f^*(K_{\mathbf{P}(X)^{\sharp}} + \Delta^{\sharp})$ 
with some divisor $\Delta_Y \ge 0$.  

The finite covering $f: Y \to \mathbf{P}(X)^{\sharp}$ can be compactified to a 
finite covering $\bar{f}: \bar{Y} \to \mathbf{P}(X)$. Let $\Delta_{\bar{Y}}$ be the closure of 
$\Delta_Y$ in $\bar{Y}$. 
Since $\mathbf{P}(X) - \mathbf{P}(X)^{\sharp}$ has codimension at least 2, one can write  $K_{\bar{Y}} + \Delta_{\bar{Y}} = \bar{f}^*(K_{\mathbf {P}(X)} + \Delta)$. 
Then, by [Xu, Proposition 1], the degree of such $\bar{f}$ is bounded by a constant 
only depending on $(\mathbf{P}(X), \Delta)$. Thus $\mathrm{deg}(f)$ is bounded above. 
Q.E.D. \vspace{0.2cm}

{\bf  3. Complex analytic orbifolds} 

We can define a complex analytic orbifold structure just by replacing the algebraic orbifold 
charts in \S 1 with the complex analytic orbifold charts (cf. [T, Chapter 13]).  
Let $U$ be a sufficiently small open neighborhood of $0 \in \mathbf{C}^m$ where 
a finite group $\Gamma$ acts on $U$ fixing the origin.   
Then a holomorphic map $\pi : U \to \mathbf{P}(X)^{\sharp}$ is called an orbifold chart 
if it is factorized as $U \to U/\Gamma \subset \mathbf{P}(X)^{\sharp}$. 
We say that two charts $\pi: U \to \mathbf{P}(X)^{\sharp}$ and $\pi': U' \to 
\mathbf{P}(X)^{\sharp}$ are admissible if $(U \times_{\mathbf{P}(X)^{\sharp}} U')^n 
\to U$ and $(U \times_{\mathbf{P}(X)^{\sharp}} U')^n \to U'$ are both \'{e}tale maps. 
Orbifold covering charts of $\mathbf{P}(X)^{\sharp}$ is a collection 
$\{\pi: U \to \mathbf{P}(X)^{\sharp}\}$ of mutually admissible charts such that $\cup \mathrm{Im}(\pi) = \mathbf{P}(X)^{\sharp}$. An orbifold structure 
on $\mathbf{P}(X)^{\sharp}$ is nothing but an equivalence class of such collections.
\vspace{0.2cm}

It is easily checked that a smooth algebraic orbifold structure on an algebraic variety $Y$  
naturally induces a complex analytic orbifold structure on $Y$. 
Conversely, if an algebraic variety $Y$ has a smooth complex analytic orbifold structure, 
then $Y$ admits a smooth algebraic orbifold structure. In fact, let $\nu: V \to V/G \subset 
Y$ be a complex analytic orbifold chart, where $V$ is an open neighborhood of 
$0 \in \mathbf{C}^m$ and $G$ fixes the origin. We put $y := \nu (0)$. 
By the local linearization of a finite group action 
[C, p.97] we may assume that $G$-action on $V$ is induced from linear  transformations of $\mathbf{C}^m$.  By Artin's approximation theorem [Ar, Corollary (2.6)] 
we may take a common \'{e}tale neighborhood $w \in W$ of $y \in Y$ and $\bar{0} \in \mathbf{C}^m/G$: 
$$ Y \leftarrow W \rightarrow \mathbf{C}^m/G.$$ 
Take the connected component $W'$ of $W \times_{\mathbf{C}^m/G}\mathbf{C}^m$ 
containing the point $(w, 0)$. Then one can write $W = W'/G'$ with a suitable subgroup $G'$ of 
$G$ and the composition map $W' \to W \to Y$ gives a smooth algebraic orbifold chart.   
\vspace{0.2cm}
 
Let $Y$ be a connected complex analytic space with an orbifold structure. Then a {\em covering 
map} $f: Y^{orb} \to \mathbf{P}(X)^{\sharp, orb}$ is a holomorphic map $f: Y \to \mathbf{P}(X)^{\sharp}$ of the underlying spaces such that 
\vspace{0.15cm}

(i) for any point $x \in \mathbf{P}(X)^{\sharp}$, 
there exists an admissible orbifold chart $\pi: U \to \mathbf{P}(X)^{\sharp}$ with 
$x \in \mathrm{Im}(\pi)$ and each connected component $V_i$ of $f^{-1}(\pi (U))$ can be 
written as $U/\Gamma_i$ where $\Gamma_i$ is some subgroup of $\Gamma$, 
\vspace{0.15cm}

(ii) the map $U \to U/\Gamma_i \cong V_i \subset Y$ is an admissible orbifold chart of 
$Y^{orb}$.    
\vspace{0.15cm}

Let $f: Y^{orb} \to \mathbf{P}(X)^{\sharp, orb}$ be an \'{e}tale covering of algebraic 
orbifolds. Then it induces a covering map of complex analytic orbifolds. 
In fact, let $y \in Y$ and $x := f(y) \in \mathbf{P}(X)^{\sharp}$. Choose algebraic orbifold charts $\mu: V \to Y$ and $\pi: U \to \mathbf{P}(X)^{\sharp}$ so that their images contain 
$y$ and $x$ respectively. We choose points $v \in V$ and $u \in U$ so that 
$\nu (v) = y$ and $\pi (u) = x$. We have a diagram of 
\'{e}tale maps $$ V \leftarrow (V \times_{\mathbf{P}(X)^{\sharp}}U)^n \to U,$$ 
and it induces an isomorphism of complex analytic germs $(V,v) \to (U,u)$. 
We have a commutative diagram 
\begin{equation} 
\begin{CD} 
(V,v) @>>> (U,u) \\ 
@VVV @VVV \\ 
(Y,y) @>>> (\mathbf{P}(X)^{\sharp},x).      
\end{CD} 
\end{equation}
By the assumption $(\mathbf{P}(X)^{\sharp}, x) \cong 
(U/G, \bar{u})$ with a finite group $G$. By the commutative diagram 
$(Y,y)$ can be also written as $(U/G', \bar{u})$ with a subgroup $G'$ of $G$. 
This shows that $f$ is a covering map of complex analytic orbifolds.  

Conversely, if $Y^{orb} \to \mathbf{P}(X)^{\sharp, orb}$ is a covering of 
complex analytic orbifolds with finite degree, then $Y$ is an algebraic variety with  
a smooth algebraic orbifold structure and the map $Y^{orb} \to \mathbf{P}^{\sharp, orb}$ 
is an \'{e}tale covering of algebraic orbifolds.  

Notice however that covering maps of  complex analytic orbifolds generally have infinite degrees. 
 
As in \S 2 take a point $x \in \mathbf{P}(X)^{\sharp} - \bar{D}^{\sharp}$ and 
consider all pairs $(Y^{orb}, y) $ of coverings $Y^{orb} \to \mathbf{P}(X)^{\sharp}$ 
and $y \in Y$ lying over $x$.  If $(Y^{orb}, y)$ and $(Y'^{orb}, y')$ are among them, then 
we can take a unique irreducible component $Z$ of $(Y \times_{\mathbf{P}(X)^{orb}}Y')^n$ 
passing through the point $z: = (y, y')$. Moreover there exists an orbifold structure on 
$Z$ such that the induced map $Z^{orb} \to \mathbf{P}(X)^{\sharp, orb}$ is a 
covering map of orbifolds. Such constructions enable us to take the inverse limit 
$(Y^{*,orb}, y^*)$ of  the inductive system $\{(Y^{orb}, y)\}$. 
Thurston [T, 13.2.4] has defined the orbifold fundamental group 
$\pi_1^{orb}(\mathbf{P}(X)^{\sharp, orb}, x)$ as the deck transformation 
group of $Y^* \to \mathbf{P}(X)^{\sharp}$. 
\vspace{0.2cm}

{\bf 4. Proof of Main Theorem} 

We fix a point $x \in \mathbf{P}(X)^{\sharp} - \bar{D}$ and 
$x^{\sharp} \in X^{\sharp}$ with $p^{\sharp}(x^{\sharp}) = x$.    
\vspace{0.2cm}

{\bf Lemma}. {\em There exists an exact sequence} 
$$ \pi_1(\mathbf{C}^*, x^{\sharp}) \to \pi_1(X^{\sharp}, x^{\sharp}) 
\to \pi_1^{orb}(\mathbf{P}(X)^{\sharp, orb}, x) \to 1.$$  
  
{\em Proof}. 
Define $$\mathbf{P}(X)^{\sharp}_0 := \mathbf{P}(X)^{\sharp} -  
\cup_{\alpha; m_{\alpha} > 1}\bar{D}_{\alpha}^{\sharp}$$ and 
$$X^{\sharp}_0 := (p^{\sharp})^{-1}(\mathbf{P}(X)^{\sharp}_0).$$ 
Since $X^{\sharp}_0 \to \mathbf{P}(X)^{\sharp}_0$ is a $\mathbf{C}^*$-bundle, 
we have an exact sequence 
$$ \pi_1(\mathbf{C}^*, x^{\sharp}) \to \pi_1(X^{\sharp}_0, x^{\sharp}) 
\to \pi_1(\mathbf{P}(X)^{\sharp}_0, x) \to 1,$$ 
where $\mathbf{C}^*$ is regarded as a fibre $(p^{\sharp})^{-1}(x)$. 
Put 
$$C := \mathrm{Coker}[\pi_1(\mathbf{C}^*, x^{\sharp}) \to \pi_1(X^{\sharp}, x^{\sharp})].$$ 
We want to prove that $C = \pi_1^{orb}(\mathbf{P}(X)^{\sharp, orb}, x^{\sharp})$. 
The kernel $N$ of the natural map $\pi_1(X^{\sharp}_0, x^{\sharp}) \to \pi_1(X^{\sharp}, x^{\sharp})$ is described as follows. 
Put $E_{\alpha} := (p^{\sharp})^{-1}(\bar{D}_{\alpha})^{\sharp}$ for $\alpha$ with $m_{\alpha} 
> 1$.  
Let $\beta'_{\alpha}$ be a small circle in $X^{\sharp}_0$ around a point of $E_{\alpha}$. 
Take a point $q_{\alpha} \in \beta'_{\alpha}$ and choose a path 
$t_{\alpha}$ in $X^{\sharp}_0$ connecting $x^{\sharp}$ and $q_{\alpha}$. 
We define a loop  $\beta_{\alpha}$ starting from $x^{\sharp}$ as 
$\beta_{\alpha} := t_{\alpha}^{-1} \circ \beta'_{\alpha} \circ t_{\alpha}$. 
Then $N$ is the smallest normal subgroup of $\pi_1(X^{\sharp}_0, x^{\sharp})$ 
containing the elements $[\beta_{\alpha}]$. 
 
Summing up these facts, one gets an exact commutative diagram   

\begin{equation} 
\begin{CD}
@.  1 @. 1 @. @. \\ 
@.  @VVV @VVV @. \\ 
@.  N  @>>> p^{\sharp}_*(N) @>>> 1 \\ 
@.  @VVV  @VVV  @. \\ 
\pi_1(\mathbf{C}^*, x^{\sharp}) @>>> \pi_1(X^{\sharp}_0, x^{\sharp}) 
@>{p^{\sharp}_*}>> \pi_1(\mathbf{P}(X)^{\sharp}_0, x) @>>> 1 \\ 
@V{id}VV @VVV @VVV @. \\ 
\pi_1(\mathbf{C}^*, x^{\sharp}) @>>> \pi_1(X^{\sharp}, x^{\sharp}) @>>> 
C @>>> 1 \\ 
@.  @VVV  @VVV  @. \\
@.  1 @. 1 @. @.
\end{CD} 
\end{equation} 

We next consider $\pi_1^{orb}(\mathbf{P}(X)^{\sharp, orb}, x)$. 
For each $\alpha$ with $m_{\alpha} > 1$, let $\gamma'_{\alpha}$ 
be a small circle in $\mathbf{P}(X)^{\sharp}_0$ around a 
point of $\bar{D}_{\alpha}^{\sharp}$. Take a point $p_{\alpha} \in \gamma_{\alpha}$ 
and choose a path $s_{\alpha}$ in $\mathbf{P}(X)^{\sharp}_0$ connecting 
$x$ and $p_{\alpha}$. We define a loop $\gamma_{\alpha}$ starting from $x$  as 
$\gamma_{\alpha} := s_{\alpha}^{-1} \circ \gamma'_{\alpha} \circ s_{\alpha}$.     
Let $M$ be the smallest normal subgroup of $\pi_1(\mathbf{P}(X)^{\sharp}_0, x)$ 
containing the elements $[\gamma_{\alpha}^{m_{\alpha}}]$.  
Then we have 
$$ \pi_1^{orb}(\mathbf{P}(X)^{\sharp, orb}, x) \cong  
\pi_1(\mathbf{P}(X)^{\sharp}_0, x)/M.$$ 
Since $p^{\sharp}_*([\beta_{\alpha}]) = [\gamma_{\alpha}^{m_{\alpha}}]$, we 
see that $M = p^{\sharp}_*(N)$. This implies that 
$C = \pi_1^{orb}(\mathbf{P}(X)^{\sharp, orb}, x^{\sharp})$.  Q.E.D.  
\vspace{0.2cm}
  
By Theorem and Lemma in section 2, there is a finite  \'{e}tale Galois covering 
$(Y^{orb}, y) \to (\mathbf{P}(X)^{\sharp}, x)$ such that $\hat{\pi}_1^{orb}(Y^{orb}, y) 
= \{1\}$. Put $Z := (X^{\sharp} \times_{\mathbf{P}(X)^{\sharp}} Y)^n$ and choose a point $z \in Z$ lying 
over $x^{\sharp}$ and $y$.  Then we have an exact sequence 
$$ \pi_1(\mathbf{C}^*, z) 
\to \pi_1(Z, z) \to \pi_1^{orb}(Y^{orb}, y) 
\to 1.$$ 
To prove that $\hat{\pi}_1(X^{\sharp}, x^{\sharp})$ is finite, it is enough to 
show that $\hat{\pi}_1(Z, z)$ is finite because $Z \to X^{\sharp}$ is a finite \'{e}tale 
covering in the usual sense.        

Assume that for an arbitrary positive integer $m$, there is a finite group $\Gamma$ 
with $\vert \Gamma \vert \geq m$ such that there is a surjection $\pi_1(Z, z) 
\to \Gamma$.  Put $$K := \mathrm{Ker}[\pi_1(\mathbf{C}^*, z) \to \pi_1(Z, z)].$$ 
Let $\Gamma_K$ be the 
image of the composition map $K \to \pi_1(Z, z) \to \Gamma$. 
Then the surjection $\pi_1(Z, z) \to \Gamma$ induces a surjection 
$\pi_1^{orb}(Y^{orb}, y) \to \Gamma/\Gamma_K$. But, since 
$\hat{\pi}_1^{orb}(Y^{orb}, y)$ is trivial, 
$\Gamma/\Gamma_K = 1$. Therefore the composition map 
$K \to \pi_1(Z, z) \to \Gamma$ is a surjection. 
If $K$ is a finite group, then this contradicts the assumption that $m$ can be 
arbitrary large. Hence the map 
$\pi_1(\mathbf{C}^*, z) \to \pi_1(Z, z)$ is an injection and 
$K = \mathbf{Z}$. Since $\Gamma$ is a quotient group of $\mathbf{Z}$ we have   
$\Gamma = \mathbf{Z}/l\mathbf{Z}$ for some $l \geq m$.  

As explained in Introduction, this leads to a contradiction when $Z$ is a $\mathbf{C}^*$-bundle over Y.  In a general case $Z$ is not a $\mathbf{C}^*$-bundle over $Y$, but it is an orbifold 
$\mathbf{C}^*$-bundle over $Y^{orb}$. Thus we can apply a similar argument to the orbifold 
$\mathbf{C}^*$-bundle $Z$ to get a contradiction:
    
Let us consider the finite \'{e}tale covering $Z' \to Z$ determined by the surjection 
$\pi_1(Z, z) \to \mathbf{Z}/l\mathbf{Z}$. By the definition this covering induces a 
cyclic covering $\mathbf{C}^* \to \mathbf{C}^*$ of degree $l$ for each general 
fibre of $Z \to Y$.  
Notice that $Z \to Y^{orb}$ is  
an orbifold  $\mathbf{C}^*$ bundle. Let  
$\mathcal{L}$ be the associated orbifold line bundle on $Y^{orb}$. Then $Z$ can 
be obtained from $\mathcal{L}$ by removing the zero section. The finite \'{e}tale covering 
map $Z' \to Z$ induces a cyclic covering $\mathcal{L}' \to \mathcal{L}$ of degree $l$ from an orbifold 
line bundle $\mathcal{L}'$ on $Y^{orb}$ branched along the zero section and $Z'$ is obtained from $\mathcal{L}'$ by removing 
the zero section . This fact, in particular, implies that $[\mathcal{L}] \in \mathrm{Pic}(Y^{orb})$ 
is divisible by $l$. Note that $m$ can be arbitrary large; but this is impossible because 
$\mathcal{L}^{-1}$ is an ample orbifold line bundle\footnote{As $\mathbf{P}(X)^{\sharp}$is 
obtained from a projective variety $\mathbf{P}(X)$ by removing a closed subset of codimension $\geq 2$, there is a complete curve inside $\mathbf{P}(X)^{\sharp}$; hence there is a complete curve $C$ inside $Y$. There is a positive integer $N$ such that $(\mathcal{M}. C) \in 1/N\cdot \mathbf{Z}$ for any orbifold line bundle $\mathcal{M}$ on $Y^{orb}$. For the 
orbifold line bundle $\mathcal{L}$ one can write $(\mathcal{L}. C)  = - k/N$ with a positive integer $k$. Then $[\mathcal{L}]$ is divisible at most by $k$.}. 
Therefore $\hat{\pi}_1(Z, z)$ is finite and so is $\hat{\pi}_1(X^{\sharp}, x^{\sharp})$.  
Since the natural map $\hat{\pi}_1(X^{\sharp}, x^{\sharp}) \to \hat{\pi}_1(X_{reg}, x^{\sharp})$ 
is surjective, we conclude that $\hat{\pi}_1(X_{reg}, x^{\sharp})$ is finite. 
\vspace{0.2cm}

{\bf Remark}. If $\pi_1^{orb}(\mathbf{P}(X)^{\sharp})$ is finite, then the argument above shows that 
$\pi_1(X_{reg})$ is finite.

\begin{center}
{\bf References} 
\end{center}

[Ar] Artin, M.: Algebraic approximation of structures over complete 
local rings, I.H.E.S . Publ. Math. {\bf 36} (1969) 23-58

[Be] Beauville, A.: Symplectic singularities, Invent. Math. {\bf 139} (2000),  
541-549

[C] Cartan, H.: Quotient d'un espace analytique par un groupe d'automorphismes, 
Algebraic geometry and Topology (A symposiumu in honor of S. Lefschetz), 
90-102 (1957) 
Princeton Univ. Press  

[C-M] Collingwood, D., McGovern, W.: Nilpotent orbits in 
semi-simple Lie algebra: an introduction, van Nostrand Reinhold Math. Series. 
Van Nostrand-Reinhold, New York (1993)   

[Fu] Fujino, O.: Introduction to the log minimal model 
program for log canonical pairs, arXiv: 0907.1506 

[K-M] Koll\'{a}r, J., Mori, S.: Birational geometry of algebraic varieties,  
Cambridge Tracts in Mathematics, {\bf 134}. Cambridge University Press, Cambridge, 1998 

[Mum] Mumford, D.: Towards an enumerative geometry of the moduli 
space of curves, in Arithmetic and Geometry, edited by M. Artin and J. Tate, 
Birkhauser-Boston, 1983, 271-326 
 
[Na] Namikawa, Y.: Equivalence of symplectic singularities, to appear in Kyoto J. Math. 
{\bf 53}, no. 2, 2013 (Memorial issue of Prof. Maruyama) 

[Sl] Slodowy, P.: Simple singularities and simple algebraic groups, Lecture Notes in 
Math. {\bf 815}. Springer, New York (1980)  

[T] Thurston, W.: The geometry and topology of 3-manifolds, 
Mimeographed Notes, Princeton Univ. Press, 1979

[Xu] Xu, C.: Finiteness of algebraic fundamental groups, 
arXiv: 1210.5564

\vspace{0.2cm}

\begin{center}
Department of Mathematics, Faculty of Science, Kyoto University, Japan 

namikawa@math.kyoto-u.ac.jp
\end{center}

\end{document}